\documentclass[12pt]{article} 
\def\R{\hbox{{\rm I}\kern-0.2em{\rm R}\kern0.2em}}%mathematical R for reals
\def\D{\hbox{{\rm I}\kern-0.2em{\rm D}\kern0.2em}}

\def\be{\begin{equation}}
\def\ee{\end{equation}}

\def\({\left(}
\def\){\right)}
\def\[{\left[}
\def\]{\right]}
\def\bc{\begin{center}}
\def\ec{\end{center}}

\begin{document}

{\large \bf Linearizability of Systems of Ordinary Differential
Equations Obtained by Complex Symmetry Analysis}

\bc{M. Safdar$^{a}$, Asghar Qadir$^{a}$, S. Ali$^{b,c}$\\
$^{a}$Center For Advanced Mathematics and Physics, National University of Sciences and Technology, Campus H-12, 44000, Islamabad, Pakistan\\
$^{b}$School of Electrical Engineering and Computer Science, National University of Sciences and Technology, Campus H-12, 44000, Islamabad, Pakistan\\
$^{c}$Present Address: Department of Mathematics, Brock University, L2S3A1 Canada \\
safdar.camp@gmail.com, aqadirmath@yahoo.com, sajid\_ali@mail.com}\ec

{\bf Abstract}. Five equivalence classes had been found for systems
of two second-order ordinary differential equations, transformable
to linear equations (linearizable systems) by a change of variables
\cite{sb}. An ``optimal (or simplest) canonical form" of linear
systems had been established to obtain the symmetry structure,
namely with 5, 6, 7, 8 and 15 dimensional Lie algebras. For those
systems that arise from a scalar complex second-order ordinary
differential equation, treated as a pair of real ordinary
differential equations, a ``reduced optimal canonical form" is
obtained. This form yields three of the five equivalence classes of
linearizable systems of two dimensions. We show that there exist
$6$, $7$ and $15$-dimensional algebras for these systems and
illustrate our results with examples.\\


\section{Introduction}
Lie used algebraic symmetry properties of differential equations to
extract their solutions \cite{lie1, lie, lie2, lie3}. One method
developed was to transform the equation to linear form by changing
the dependent and independent variables invertibly. Such
transformations are called \emph{point transformations} and the
transformed equations are said to be \emph{linearized}. Equations
that can be so transformed are said to be {\it linearizable}. Lie
proved that the necessary and sufficient condition for a scalar
nonlinear ordinary differential equation (ODE) to be linearizable is
that it must have eight Lie point symmetries. He exploited the fact
that all scalar linear second-order ODEs are equivalent under point
transformations \cite{sur}, i.e. every linearizable scalar
second-order ODE is reducible to the free particle equation. While
the situation is not so simple for scalar linear ODEs of order
$n\geq3$, it was proved that there are three equivalence classes
with $n+1$, $n+2$ or $n+4$ infinitesimal symmetry generators
\cite{lea}.

For linearization of systems of two nonlinear ODEs, we will first
consider the equivalence of the corresponding linear systems under
point transformations. Nonlinear systems of two second-order ODEs
that are linearizable to systems of ODEs with constant coefficients,
were proved to have three equivalence classes \cite{gor}. They have
$7$, $8$ or $15$-dimensional Lie algebras. This result was extended
to those nonlinear systems which are equivalent to linear systems of
ODEs with constant or variable coefficients \cite{sb}. They obtained
an ``optimal" canonical form of the linear systems involving three
parameters, whose specific choices yielded five equivalence classes,
namely with $5$, $6$, $7$, $8$ or $15$-dimensional Lie algebras.

Geometric methods were developed to transform nonlinear systems of
second-order ODEs \cite{aa,mq1,mq2} to a system of the free particle
equations by treating them as geodesic equations and then projecting
those equations down from an $m\times m$ system to an $(m-1)\times
(m-1)$ system. In this process the originally homogeneous
quadratically semi-linear system in $m$ dimensions generically
becomes a non-homogeneous, cubically semi-linear system in $(m-1)$
dimensions. When used for $m=2$ the Lie conditions for the scalar
ODE are recovered precisely. The criterion for linearizability is
simply that the manifold for the (projected) geodesic equations be
flat. The symmetry algebra in this case is $sl(n+2,\R)$ and hence
the number of generators is $n^2+4n+3$. Thus for a system of two
equations to be linearizable by this method it must have 15
generators.

A scalar complex ODE involves two real functions of two real
variables, yielding a system of two partial differential equations
(PDEs) \cite{saj, saj1}. By restricting the  independent variable to
be real we obtain a system of ODEs. Complex symmetry analysis (CSA)
provides the symmetry algebra for systems of two ODEs with the help
of the symmetry generators of the corresponding complex ODE. This is
not a simple matter of doubling the generators for the scalar
complex ODE. The inequivalence of these systems with the above
mentioned systems obtained earlier (by geometric means) \cite{mq2},
has been proved \cite{saf2}. \emph{Thus their symmetry structures
are not the same}. We prove that a general two-dimensional system of
second-order ODEs corresponds to a scalar complex second-order ODE
if the coefficients of the system satisfy Cauchy-Riemann equations
(CR-equations). We provide the full symmetry algebra for the systems
of ODEs that correspond to linearizable scalar complex ODEs. For
this purpose we derive a \emph{reduced optimal canonical form} for
linear systems obtainable from a complex linear equation. We prove
that this form provides three equivalence classes of linearizable
systems of two second-order ODEs, while there exist five
linearizable classes \cite{sb} by real symmetry analysis. This
difference arises due to the fact that in CSA we invoke equivalence
of {\it scalar} second-order ODEs to obtain the reduced optimal
form, while in real symmetry analysis equivalence of linear {\it
systems} of two ODEs was used to derive their optimal form. The
nonlinear systems transformable to one of the three equivalence
classes we provide here, are characterized by complex
transformations of the form
\begin{eqnarray*}
T:(x,u(x))\rightarrow (\chi(x),U(x,u)).
\end{eqnarray*}
Indeed, these complex transformations generate these linearizable
classes of two dimensional systems. Note that not all the complex
linearizing transformations for scalar complex equations provide the
corresponding real transformations for systems.

The plan of the paper is as follows. In the next section we present
the preliminaries for determining the symmetry structures. The third
section deals with the conditions derived for systems that can be
obtained by CSA. In section four we obtain the reduced optimal
canonical form for systems associated with complex linear ODEs. The
theory developed to classify linearizable systems of ODEs
transformable to this reduced optimal form is given in the fifth
section. Applications of the theory are given in the next section.
The last section summarizes and discusses the work.

%%%%%%%%%%%%%%%%%%%%%%%%%%%%%%%%%%%%%%%%%%%%%%%%%%%%%%%%%%%%%%%%%%%%%%%%%%%%%%%%%%%%%%%%%%%%%%%%%%%%%%%%%%%%%%%%%%%%%%%
%%%%%%%%%%%%%%%%%%%%%%%%%%%%%%%%%%%%%%%%%%%%%%%%%%%%%%%%%%%%%%%%%%%%%%%%%%%%%%%%%%%%%%%%%%%%%%%%%%%%%%%%%%%%%%%%%%%%%%%

\section{Preliminaries}

The simplest form of a second-order equation has the
maximal-dimensional algebra, $sl(3,\R)$. To discuss the equivalence
of systems of two linear second-order ODEs, we need to use the
following result for the equivalence of a general system of $n$
linear homogeneous second-order ODEs with $2n^{2}+n$ arbitrary
coefficients and some canonical forms that have fewer arbitrary
coefficients \cite{wf}. Any system of $n$ second-order
non-homogeneous linear ODEs
\begin{eqnarray}
\ddot{\textbf{u}}=\textbf{A} \dot{\textbf{u}}+\textbf{B}
\textbf{u}+\textbf{c},\label{1}
\end{eqnarray}
can be mapped invertibly to one of the following forms
\begin{equation}
\ddot{\textbf{v}}=\textbf{C} \dot{\textbf{v}},\label{2}
\end{equation}
\begin{equation}
\ddot{\textbf{w}}=\textbf{D} \textbf{w},\label{3}
\end{equation}
where $\textbf{A}$, $\textbf{B}$, $\textbf{C}$, $\textbf{D}$ are
$n\times n$ matrix functions, $\textbf{u}$, $\textbf{v}$,
$\textbf{w}$, $\textbf{c}$ are vector functions and dot represents
differentiation relative to the independent variable $t$. For a
system of two second-order ODEs ($n=2$) there are a total of $10$
coefficients for the system represented by equation ($\ref{1}$). It
is reducible to the first and second canonical forms, ($\ref{2}$)
and ($\ref{3}$) respectively. Thus a system with $4$ arbitrary
coefficients of the form
\begin{eqnarray}
\ddot{w_{1}}=d_{11}(t)w_{1}+d_{12}(t)w_{2},\nonumber\\
\ddot{w_{2}}=d_{21}(t)w_{1}+d_{22}(t)w_{2},\label{4}
\end{eqnarray}
can be obtained by using the equivalence of ($\ref{1}$) and the
counterpart of the Laguerre-Forsyth second canonical form
($\ref{3}$). This result demonstrates the equivalence of systems of
two ODEs having $10$ and $4$ arbitrary coefficients respectively.
The number of arbitrary coefficients can be further reduced to three
by the change of variables \cite{sb}
\begin{eqnarray}
\tilde{y}=w_{1}/\rho(t),\hspace{2mm}\tilde{z}=w_{2}/\rho(t),
\hspace{2mm}x=\int^{t}\rho^{-2}(s)ds,\label{5}
\end{eqnarray}
where $\rho$ satisfies
\begin{equation}
\rho^{\prime\prime}-\frac{d_{11}+d_{22}}{2}\rho=0,\label{6}
\end{equation}
to the linear system
\begin{eqnarray}
\tilde{y}^{\prime\prime}=\tilde{d}_{11}(x)\tilde{y}+\tilde{d}_{12}(x)\tilde{z},\nonumber\\
\tilde{z}^{\prime\prime}=\tilde{d}_{21}(x)\tilde{y}-\tilde{d}_{11}(x)\tilde{z},\label{7}
\end{eqnarray}
where
\begin{eqnarray}
\tilde{d}_{11}=\frac{\rho^{3}(d_{11}-d_{22})}{2},\hspace{2mm}\tilde{d}_{12}=
\rho^{3}d_{12},\hspace{2mm}\tilde{d}_{21}=\rho^{3}d_{21}.\label{8}
\end{eqnarray}
This procedure of reduction of arbitrary coefficients for
linearizable systems simplifies the classification problem
enormously. System ($\ref{7}$) is called the \emph{optimal canonical
form} for linear systems of two second-order ODEs, as it has the
fewest arbitrary coefficients, namely three.
%%%%%%%%%%%%%%%%%%%%%%%%%%%%%%%%%%%%%%%%%%%%%%%%%%%%%%%%%%%%%%%%%%%%%%%%%%%%%%%%%%%%%%%%%%%%%%%%%%%%%%%%%%%%%%%%%%%%%%%
%%%%%%%%%%%%%%%%%%%%%%%%%%%%%%%%%%%%%%%%%%%%%%%%%%%%%%%%%%%%%%%%%%%%%%%%%%%%%%%%%%%%%%%%%%%%%%%%%%%%%%%%%%%%%%%%%%%%%%%

\section{Systems of ODEs obtainable by CSA}

Following the classical Lie procedure, one uses point
transformations
\begin{equation}
X=X(x,y,z),\hspace{2mm}Y=Y(x,y,z),\hspace{2mm}Z=Z(x,y,z),\label{9}
\end{equation}
to map the general linearizable system of two second-order ODEs
\cite{jm}, which is (at most) cubically semi-linear in both the
dependent variables,
\begin{eqnarray}
y^{^{\prime \prime }}=\omega_{1}(x,y,z,y^{^{\prime }},z^{^{\prime }}),\nonumber \\
z^{^{\prime \prime }}=\omega_{2}(x,y,z,y^{^{\prime }},z^{^{\prime
}}),\label{10}
\end{eqnarray}
where prime denotes differentiation relative to $x$, to the simplest
form
\begin{equation}
Y^{\prime\prime}=0,\hspace{2mm}Z^{\prime\prime}=0,\label{11}
\end{equation}
where the prime now denotes differentiation with respect to $X$ and
the mappings ($\ref{9}$) are invertible. The derivatives transform
as
\begin{eqnarray}
Y^{\prime}=\frac{D_{x}(Y)}{D_{x}(X)}=F_{1}(x,y,z,y^{^{\prime
}},z^{^{\prime }}),\nonumber \\
Z^{\prime}=\frac{D_{x}(Z)}{D_{x}(X)}=F_{2}(x,y,z,y^{^{\prime
}},z^{^{\prime}}),\label{12}
\end{eqnarray}
and
\begin{equation}
Y^{\prime\prime}=\frac{D_{x}(F_{1})}{D_{x}(X)},~~
Z^{\prime\prime}=\frac{D_{x}(F_{2})}{D_{x}(X)},\label{13}
\end{equation}
where $D_{x}$ is the total derivative operator. This yields
\begin{equation}
\begin{tabular}{l}
$y^{^{\prime \prime }}+\alpha_{11}y^{^{\prime
}3}+\alpha_{12}y^{^{\prime }2}z^{^{\prime }}+\alpha_{13}y^{^{\prime
}}z^{^{\prime }2}+\alpha_{14}z^{^{\prime }3}+\beta_{11}y^{^{\prime
}2}+\beta_{12}y^{^{\prime
}}z^{^{\prime }}+\beta_{13}z^{^{\prime }2}$ \\
$+\gamma_{11}y^{^{\prime }}+\gamma_{12}z^{^{\prime }}+\delta_{1}=0,$ \\
\\
$z^{^{\prime \prime }}+\alpha_{21}y^{^{\prime
}3}+\alpha_{22}y^{^{\prime }2}z^{^{\prime }}+\alpha_{23}y^{^{\prime
}}z^{^{\prime }2}+\alpha_{24}z^{^{\prime }3}+\beta_{21}y^{^{\prime
}2}+\beta_{22}y^{^{\prime
}}z^{^{\prime }}+\beta_{23}z^{^{\prime }2}$ \\
$+\gamma_{21}y^{^{\prime }}+\gamma_{22}z^{^{\prime
}}+\delta_{2}=0,$\label{14}
\end{tabular}
\end{equation}
the coefficients being functions of the independent and dependent
variables. System ($\ref{14}$) is the most general candidate for two
second-order ODEs that may be linearizable. While another candidate
of linearizability of two dimensional systems obtainable from the
most general form of a complex linearizable equation
\begin{eqnarray}
u^{\prime\prime}+E_{3}(x,u)u^{\prime 3}+E_{2}(x,u)u^{\prime
2}+E_{1}(x,u)u^{\prime}+E_{0}(x,u)=0,\label{15}
\end{eqnarray}
where $u$ is a complex function of the real independent variable
$x$, is also cubically semi-linear i.e. a system of the form
\begin{equation}
\begin{tabular}{l}
$y^{^{\prime \prime }}+\bar \alpha_{11}y^{^{\prime }3}-3\bar
\alpha_{12}y^{^{\prime }2}z^{^{\prime }}-3\bar
\alpha_{11}y^{^{\prime }}z^{^{\prime }2}+\bar \alpha_{12}z^{^{\prime
}3}+\bar \beta_{11}y^{^{\prime }2}-2\bar \beta_{12}y^{^{\prime
}}z^{^{\prime }}-\bar \beta_{11}z^{^{\prime }2}$ \\
$+\bar \gamma_{11}y^{^{\prime }}-\bar \gamma_{12}z^{^{\prime }}+\bar \delta_{11}=0,$ \\
\\
$z^{^{\prime \prime }}+\bar \alpha_{12}y^{^{\prime }3}+3\bar
\alpha_{11}y^{^{\prime }2}z^{^{\prime }}-3\bar
\alpha_{12}y^{^{\prime }}z^{^{\prime }2}-\bar \alpha_{11}z^{^{\prime
}3}+\bar \beta_{12}y^{^{\prime }2}+2\bar \beta_{11}y^{^{\prime
}}z^{^{\prime }}-\bar \beta_{12}z^{^{\prime }2}$ \\
$+\bar \gamma_{12}y^{^{\prime }}+\bar \gamma_{11}z^{^{\prime }}+\bar
\delta_{12}=0,$\label{16}
\end{tabular}
\end{equation}
here the coefficients $\bar \alpha_{1i}$, $\bar \beta_{1i}$, $\bar
\gamma_{1i}$ and $\bar \delta_{1i}$ for $i=1,2$ are functions of
$x$, $y$ and $z$. Clearly, the system ($\ref{16}$) corresponds to
($\ref{15}$) if the coefficients $\bar \alpha_{1i},~\bar
\beta_{1i},~\bar \gamma_{1i}$ and $\bar \delta_{1i}$ satisfy the
CR-equations i.e.
$\alpha_{11,y}=\alpha_{12,z},~\alpha_{12,y}=-\alpha_{11,z}$ and vice
versa. It is obvious as ($\ref{15}$) generates a system by breaking
the complex coefficients $E_{j}$, for $j=0,1,2,3$ into real and
imaginary parts
\begin{eqnarray}
E_{3}=\bar \alpha_{11}+i\bar \alpha_{12},~~ E_{2}=\bar
\beta_{11}+i\bar \beta_{12},~~ E_{1}=\bar \gamma_{11}+i\bar
\gamma_{12},~~ E_{0}=\bar \delta_{11}+i\bar \delta_{12},\label{17}
\end{eqnarray}
where all the functions are analytic. Hence we can state the following theorem.\\
\newline \textbf{Theorem 1.} \textit{A general two dimensional system of second-order ODEs} ($\ref{10}$) \textit{corresponds to a complex equation}
\begin{eqnarray}
u^{\prime\prime}=\omega(x,u,u^{\prime}),
\end{eqnarray}
\textit{if and only if $\omega_{1}$ and $\omega_{2}$ satisfy the
CR-equations}
\begin{eqnarray}
\omega_{1,y}=\omega_{2,z},~~\omega_{1,z}=-\omega_{2,y},\nonumber\\
\omega_{1,y^{\prime}}=\omega_{2,z^{\prime}},~~\omega_{1,z^{\prime}}=-\omega_{2,y^{\prime}}.
\end{eqnarray}
\\
For the correspondence of both the cubic forms ($\ref{14}$) and
($\ref{16}$) of two dimensional systems we state the following
theorem.\\
\newline \textbf{Theorem 2.} \textit{A system of the
form} ($\ref{14}$) \textit{corresponds to} ($\ref{16}$) \textit{if
and only if the coefficients $\alpha_{ij}$, $\beta_{ik}$,
$\gamma_{il}$ and $\delta_{i}$ satisfy the following conditions
\begin{eqnarray}
\alpha_{11}=-\frac{1}{3}\alpha_{13}=\frac{1}{3}\alpha_{22}=-\alpha_{24},\nonumber\\
-\frac{1}{3}\alpha_{12}=\alpha_{14}=\alpha_{21}=-\frac{1}{3}\alpha_{23},\nonumber\\
\beta_{11}=\frac{1}{2}\beta_{22}=-\beta_{13},\nonumber\\
\beta_{21}=-\frac{1}{2}\beta_{12}=-\beta_{23},\nonumber\\
\gamma_{11}=\gamma_{22}=,\quad \gamma_{21}=-\gamma_{12},\label{18}
\end{eqnarray}
where $i=l=1,2$,  $j=1,...,4$ and $k=1,2,3$.}\\
\newline \textbf{Proof.} It can be trivially proved if we rewrite the above
equations as $\bar \alpha_{1i}$, $\bar \beta_{1i}$ and $\bar
\gamma_{1i}$, respectively. These coefficients correspond to complex
coefficients of ($\ref{15}$) if and only if they satisfy the
CR-equations.\\
\newline Thus Theorem (1) and (2) identify those two dimensional systems
which are obtainable from complex equations.
%%%%%%%%%%%%%%%%%%%%%%%%%%%%%%%%%%%%%%%%%%%%%%%%%%%%%%%%%%%%%%%%%%%%%%%%%%%%%%%%%%%%%%%%%%%%%%%%%%%%%%%%%%%%%%%%%%%%%%%
%%%%%%%%%%%%%%%%%%%%%%%%%%%%%%%%%%%%%%%%%%%%%%%%%%%%%%%%%%%%%%%%%%%%%%%%%%%%%%%%%%%%%%%%%%%%%%%%%%%%%%%%%%%%%%%%%%%%%%%

\section{Reduced optimal canonical forms}

The simplest forms for linear systems of two second-order ODEs
corresponding to complex scalar ODEs can be established by invoking
the equivalence of scalar second-order linear ODEs. Consider a
general linear scalar complex second-order ODE
\begin{equation}
u^{\prime\prime}=\zeta_{1}(x)u^{\prime}+\zeta_{2}(x)u+\zeta_{3}(x),\label{21}
\end{equation}
where prime denotes differentiation relative to $x$ and
$u(x)=y(x)+iz(x)$ is a complex function of the real independent
variable $x$. As all the linear scalar second-order ODEs are
equivalent, so equation ($\ref{21}$) is equivalent to the following
scalar second-order complex ODEs
\begin{equation}
u^{\prime\prime}=\zeta_{4}(x)u^{\prime},\label{22}
\end{equation}
\begin{equation}
u^{\prime\prime}=\zeta_{5}(x)u,\label{23}
\end{equation}
where all the three forms ($\ref{21}$), ($\ref{22}$) and
($\ref{23}$) are transformable to each other. Indeed these three
forms are reducible to the free particle equation. These three
complex scalar linear ODEs belong to the same equivalence class,
i.e. all have eight Lie point symmetry generators. In this paper we
prove that the systems obtainable by these forms using CSA have more
than one equivalence class. To extract systems of two linear ODEs
from ($\ref{22}$) and ($\ref{23}$) we put
$\zeta_{4}(x)=\alpha_{1}(x)+i\alpha_{2}(x)$ and
$\zeta_{5}(x)=\alpha_{3}(x)+i\alpha_{4}(x)$ to obtain two linear
forms of system of two linear second-order ODEs
\begin{eqnarray}
y^{\prime\prime}=\alpha_{1}(x)y^{\prime}-\alpha_{2}(x)z^{\prime},\nonumber\\
z^{\prime\prime}=\alpha_{2}(x)y^{\prime}+\alpha_{1}(x)z^{\prime}.\label{24}
\end{eqnarray}
and
\begin{eqnarray}
y^{\prime\prime}=\alpha_{3}(x)y-\alpha_{4}(x)z,\nonumber\\
z^{\prime\prime}=\alpha_{4}(x)y+\alpha_{3}(x)z,\label{25}
\end{eqnarray}
thus we state the following theorem.\\
\newline \textbf{Theorem 3.}
\textit{If a system of two second-order ODEs is linearizable via
invertible complex point transformations then it can be mapped to
one of the two forms} ($\ref{24}$) \textit{or} ($\ref{25}$).

Notice that here we have \emph{only two arbitrary coefficients} in
both the linear forms, while the minimum number obtained before was
three i.e. a system of the form ($\ref{7}$). The reason we can
reduce further is that we are dealing with the special classes of
linear systems of ODEs that correspond to the scalar complex ODEs.
In fact ($\ref{25}$) can be reduced further by the change of
variables
\begin{eqnarray}
Y=y/\rho(t),\hspace{2mm}Z=z/\rho(t),\hspace{2mm}
x=\int^{t}\rho^{-2}(s)ds,\label{26}
\end{eqnarray}
where $\rho$ satisfies
\begin{equation}
\rho^{\prime\prime}-\alpha_{3}\rho=0,\label{27}
\end{equation}
to
\begin{eqnarray}
Y^{\prime\prime}=-\beta(x)Z,~~~~
Z^{\prime\prime}=\beta(x)Y,\label{28}
\end{eqnarray}
where $\beta=\rho^{3}\alpha_{4}$. We state this result in the form
of a theorem.\\
\newline \textbf{Theorem 4.} \textit{Any linear
system of two second-order ODEs of the form} ($\ref{25}$)
\textit{with two arbitrary coefficients is transformable to a
simplest system of two linear ODEs} ($\ref{28}$) \textit{with one
arbitrary coefficient via real point transformations} ($\ref{26}$)
\textit{and} ($\ref{27}$).

Equation ($\ref{28}$) is the \emph{reduced optimal canonical form}
for systems associated with complex ODEs, with just one coefficient
which is an arbitrary function of $x$. The equivalence of systems
($\ref{24}$) and ($\ref{25}$) can be established via invertible
point transformations, so we state the following theorem.\\
\newline
\textbf{Theorem 5.} \textit{Two linear forms of the systems of two
second-order ODEs} ($\ref{24}$) \textit{and} ($\ref{25}$)
\textit{are equivalent via invertible point transformations}
\begin{eqnarray}
y=M_{1}(x)y_{1}-M_{2}(x)y_{2}+y^{*},\nonumber\\
z=M_{1}(x)y_{2}+M_{2}(x)y_{1}+z^{*},\label{29}
\end{eqnarray}
\textit{of the dependent variables only, where} $M_{1}(x)$,
$M_{2}(x)$ \textit{are two linearly independent solutions of}
\begin{eqnarray}
\alpha_{1}M_{1}-\alpha_{2}M_{2}=2M^{\prime}_{1},\nonumber\\
\alpha_{1}M_{1}-\alpha_{2}M_{2}=2M^{\prime}_{1},\label{30}
\end{eqnarray}
\textit{and} $y^{*}$, $z^{*}$ \textit{are the particular solutions
of} ($\ref{24}$).\\
\newline {\bf Proof.} Differentiating the
set of equations ($\ref{30}$) and using the result in the linear
form ($\ref{24}$), Routine calculations show that ($\ref{24}$) can
be mapped to ($\ref{25}$) where
\begin{eqnarray}
\alpha_{3}(x)=\frac{1}{M^{2}_{1}+M^{2}_{2}}[M_{1}(\alpha_{1}M^{\prime}_{1}-\alpha_{2}M^{\prime}_{2}-M^{\prime\prime}_{1})+M_{2}(\alpha_{1}M^{\prime}_{2}+\alpha_{2}M^{\prime}_{1}-M^{\prime\prime}_{2})],\nonumber\\
\alpha_{4}(x)=\frac{1}{M^{2}_{1}+M^{2}_{2}}[M_{1}(\alpha_{1}M^{\prime}_{2}+\alpha_{2}M^{\prime}_{1}-M^{\prime\prime}_{2})-M_{2}(\alpha_{1}M^{\prime}_{1}-\alpha_{2}M^{\prime}_{2}-M^{\prime\prime}_{1})].\nonumber\\
\label{31}
\end{eqnarray}
Thus the linear form ($\ref{24}$) is reducible to
($\ref{28}$).\\
\newline \textbf{Remark 1.} Any nonlinear system of
two second-order ODEs that is linearizable by complex methods can be
mapped invertibly to a system of the form ($\ref{28}$) with one
coefficient which is an arbitrary function of the independent
variable.
%%%%%%%%%%%%%%%%%%%%%%%%%%%%%%%%%%%%%%%%%%%%%%%%%%%%%%%%%%%%%%%%%%%%%%%%%%%%%%%%%%%%%%%%%%%%%%%%%%%%%%%%%%%%%%%%%%%%%%%
%%%%%%%%%%%%%%%%%%%%%%%%%%%%%%%%%%%%%%%%%%%%%%%%%%%%%%%%%%%%%%%%%%%%%%%%%%%%%%%%%%%%%%%%%%%%%%%%%%%%%%%%%%%%%%%%%%%%%%%
\section{Symmetry structure of linear systems obtained by CSA}

To use the reduced canonical form \cite{hs} for deriving the
symmetry structure of linearizable systems associated with the
complex scalar linearizable ODEs, we obtain a system of PDEs whose
solution provides the symmetry generators for the corresponding
linearizable systems of
two second-order ODEs.\\
\newline \textbf{Theorem 6.} \textit{Linearizable systems of two second-order ODEs reducible to the
linear form} ($\ref{28}$) \textit{via invertible complex point
transformations, have} $6$, $7$ \textit{or
$15$-dimensional Lie point symmetry algebras.}\\
\newline
\textbf{Proof.} The symmetry conditions provide the following set of
PDEs for the system ($\ref{28}$)
\begin{eqnarray}
\xi_{xx}=\xi_{xy}=\xi_{yy}=0=\eta_{1,zz}=\eta_{2,yy},\label{32}\\
\eta_{1,yy}-2\xi_{xy}=\eta_{1,yz}-\xi_{xz}=\eta_{2,yz}-\xi_{xy}=
\eta_{2,zz}-2\xi_{xz}=0,\label{33}\\
\xi_{xx}-2\eta_{1,xy}-3\beta(x)\xi_{,y}z+\beta(x)\xi_{,z}y
=\eta_{1,xz}+\beta(x)\xi_{,z}z=0,\label{34}\\
\xi_{xx}-2\eta_{2,xz}+3\beta(x)\xi_{,z}y-\beta(x)\xi_{,y}z=
\eta_{2,xy}-\beta(x)\xi_{,y}y=0,\label{35}\\
\eta_{1,xx}+\beta(x)(\eta_{1,z}y+2\xi_{,x}z-\eta_{1,y}z+\eta_{2})
+\beta^{\prime}(x)z\xi=0,\label{36}\\
\eta_{2,xx}+\beta(x)(\eta_{2,z}y-2\xi_{,x}y-\eta_{2,y}z-\eta_{1})
-\beta^{\prime}(x)y\xi=0.\label{37}
\end{eqnarray}
Equations ($\ref{34}$)-($\ref{37}$) involve an arbitrary function of
the independent variable and its first derivatives. Using equations
($\ref{32}$) and ($\ref{33}$) we have the following solution set
\begin{eqnarray}
\xi=\gamma_{1}(x)y+\gamma_{2}(x)z+\gamma_{3}(x),\nonumber\\
\eta_{1}=\gamma^{\prime}_{1}(x)y^{2}+\gamma^{\prime}_{2}(x)yz+
\gamma_{4}(x)y+\gamma_{5}(x)z+\gamma_{6}(x),\nonumber\\
\eta_{2}=\gamma^{\prime}_{1}(x)yz+\gamma^{\prime}_{2}(x)z^{2}+
\gamma_{7}(x)y+\gamma_{8}(x)z+\gamma_{9}(x).\label{38}
\end{eqnarray}
Using equations ($\ref{34}$) and ($\ref{35}$), we get
\begin{eqnarray}
\beta(x)\gamma_{1}(x)=0=\beta(x)\gamma_{2}(x).\label{39}
\end{eqnarray}
Now assuming $\beta(x)$ to be zero, non-zero constant and arbitrary function of $x$ will generate the following cases.\\
\newline\textbf{Case 1.1.} $\beta(x)=0$.\\
The set of determining equations ($\ref{32}$)-($\ref{37}$) will
reduce to a trivial system of PDEs
\begin{eqnarray}
\eta_{1,xx}=\eta_{1,xz}=\eta_{1,zz}=0,\nonumber\\
\eta_{2,xx}=\eta_{2,xy}=\eta_{2,yy}=0,\nonumber\\
2\xi_{,xy}-\eta_{1,yy}=0=2\xi_{,xz}-\eta_{2,zz},\nonumber\\
\xi_{,xz}-\eta_{1,yz}=0=\xi_{,xy}-\eta_{2,yz},\nonumber\\
\xi_{,xx}-2\eta_{1,xy}=0=\xi_{,xx}-2\eta_{2,xz},\label{40}
\end{eqnarray}
which can be extracted classically for the system of free particle
equations. Solving it we find a $15$-dimensional Lie point symmetry algebra.\\
\newline\textbf{Case 1.2.} $\beta(x)\neq0$.\\
Then ($\ref{39}$) implies $\gamma_{1}(x)=\gamma_{2}(x)=0$ and
($\ref{38}$) reduces to
\begin{eqnarray}
\xi=\gamma_{3}(x),\nonumber\\
\eta_{1}=(\frac{\gamma^{\prime}_{3}(x)}{2}+c_{3})y+c_{1}z+
\gamma_{6}(x),\nonumber\\
\eta_{2}=c_{2}y+(\frac{\gamma^{\prime}_{3}(x)}{2}+c_{4})z+
\gamma_{9}(x).\label{41}
\end{eqnarray}
Here two subcases arise.\\
\newline\textbf{Case 1.2.1.} $\beta(x)$ \emph{is a non-zero constant}.\\
As equations ($\ref{36}$) and ($\ref{37}$) involve the derivatives
of $\beta(x)$, which will now be zero, equations
($\ref{34}$)-($\ref{37}$) and ($\ref{41}$) yield a $7$-dimensional
Lie algebra. The explicit expressions of the symmetry generators
involve trigonometric functions. But for a simple demonstration of
the algorithm consider $\beta(x)=1$. The solution of the set of the
determining equations is
\begin{equation}
\xi=C_{1},\nonumber
\end{equation}
and
\begin{eqnarray}
\eta_{1}=C_{2}y+[-C_{4}e^{x/\sqrt{2}}-C_{3}e^{-x/\sqrt{2}}]
\sin({x/\sqrt{2}})+C_{6}e^{x/\sqrt{2}}\cos({x/\sqrt{2}})+
\nonumber\\
C_{5}e^{-x/\sqrt{2}}\cos({x/\sqrt{2}})+C_{7}z,\nonumber\\
\eta_{2}=[-C_{6}e^{x/\sqrt{2}}+C_{5}e^{-x/\sqrt{2}}]
\sin({x/\sqrt{2}})-C_{4}e^{x/\sqrt{2}}\cos({x/\sqrt{2}})-C_{2}z+
\nonumber\\
C_{3}e^{-x/\sqrt{2}}\cos({x/\sqrt{2}})+C_{7}y.\label{42}
\end{eqnarray}
This yields a $7$-dimensional symmetry algebra.\\
\newline\textbf{Case 1.2.2.1.} $\beta(x)=x^{-2}, x^{-4}$ or $(x+1)^{-4}$.\\
Equations ($\ref{34}$)-($\ref{37}$) and ($\ref{41}$) yield a
$7$-dimensional Lie algebra. Thus the $7$-dimensional algebras can
be related with systems which have variable coefficients in their
linear forms, apart from the linear forms with constant coefficients.\\
\newline\textbf{Case 1.2.2.2.} $\beta(x)=x^{-1}, x^{2}$, $x^{2}\pm C_{0}$ or $e^{x}$.\\
Using equations ($\ref{34}$)-($\ref{37}$) and ($\ref{41}$), we
arrive at a $6$-dimensional Lie point symmetry algebra. The explicit
expressions involve special functions, e.g for $\beta(x)=x^{-1}$,
$x^{2}$, $x^{2}\pm C_{0}$ we get Bessel functions. Similarly for
$\beta(x)=e^{x}$ there are six symmetries, including the generators
$y\partial_{y}-e^{x}z\partial_{z}$,
$z\partial_{z}+e^{x}y\partial_{y}$. The remaining four generators
come from the solution of an ODE of order four.\\

Thus there is only a $6$, $7$ or $15$-dimensional algebra for
linearizable systems of two second-order ODEs transformable to
($\ref{28}$) via invertible complex point transformations. We are
not investigating the remaining two linear forms ($\ref{24}$) and
($\ref{25}$), because these are transformable to system ($\ref{28}$)
i.e. all these forms have the same symmetry structures. The linear
forms providing $6$ or $7$-dimensional algebras here are obtainable
by linear forms extractable from ($\ref{7}$), with a $6$ or $7$
dimensional algebra respectively. Consider ($\ref{7}$) with all the
coefficients to be non-zero constants i.e.
$\tilde{d}_{11}(x)=a_{0}$, $\tilde{d}_{12}(x)=b_{0}$ and
$\tilde{d}_{21}(x)=c_{0}$, where
\begin{eqnarray}
a_{0}^{2}+b_{0}c_{0}\neq 0.\label{43}
\end{eqnarray}
This system provides seven symmetry generators. The linear form
($\ref{28}$) also provides a $7$-dimensional algebra with constant
coefficients satisfying ($\ref{43}$), while the $8$-dimensional
symmetry algebra was extracted \cite{sb} by assuming
\begin{eqnarray}
a_{0}^{2}+b_{0}c_{0}=0.\label{44}
\end{eqnarray}
Such linear forms cannot be obtained from ($\ref{28}$). These two
examples explain why a $7$-dimensional algebra can be obtained from
($\ref{28}$), but a linear form with an $8$-dimensional algebra is
not obtainable from
it.\\
\newline To prove these observations consider arbitrary point
transformations of the form
\begin{eqnarray}
\tilde{y}=a(x)y+b(x)z,~~~\tilde{z}=c(x)y+d(x)z.\label{45}
\end{eqnarray}
\newline\textbf{Case a.} If $a(x)=a_{0}$, $b(x)=b_{0}$, $c(x)=c_{0}$ and
$d(x)=d_{0}$ are constants then ($\ref{45}$) implies
\begin{eqnarray}
\tilde{y}^{\prime\prime}=a_{0}y^{\prime\prime}+b_{0}z^{\prime\prime},\nonumber\\
\tilde{z}^{\prime\prime}=c_{0}y^{\prime\prime}+d_{0}z^{\prime\prime}.\label{46}
\end{eqnarray}
Using ($\ref{7}$) and ($\ref{25}$) in the above equation we find
\begin{eqnarray}
(a_{0}d_{0}-b_{0}c_{0})y^{\prime\prime}=((a_{0}d_{0}+b_{0}c_{0})\tilde{d}_{11}(x)+c_{0}d_{0}\tilde{d}_{12}(x)-a_{0}b_{0}\tilde{d}_{21}(x))y+\nonumber\\
(2b_{0}d_{0}\tilde{d}_{11}(x)+d^{2}_{0}\tilde{d}_{12}(x)-b^{2}_{0}\tilde{d}_{21}(x))z,\nonumber\\
(a_{0}d_{0}-b_{0}c_{0})z^{\prime\prime}=((a_{0}d_{0}+b_{0}c_{0})\tilde{d}_{11}(x)+c_{0}d_{0}\tilde{d}_{12}(x)-a_{0}b_{0}\tilde{d}_{21}(x))z+\nonumber \\
(2a_{0}c_{0}\tilde{d}_{11}(x)+c^{2}_{0}\tilde{d}_{12}(x)-a^{2}_{0}\tilde{d}_{21}(x))y,\label{47}
\end{eqnarray}
where $a_{0}d_{0}-b_{0}c_{0}\neq0$. Using ($\ref{25}$), ($\ref{47}$)
and the linear independence of the $\tilde{d}$'s, gives
\begin{eqnarray}
a_{0}b_{0}=c_{0}d_{0}=0,\nonumber\\
a^{2}_{0}-b^{2}_{0}=c^{2}_{0}-d^{2}_{0}=0,\nonumber\\
a_{0}d_{0}+b_{0}c_{0}=a_{0}c_{0}-b_{0}d_{0}=0,\label{48}
\end{eqnarray}
which has a solution $a_{0}=b_{0}=c_{0}=d_{0}=0$, which is
inconsistent with ($\ref{47}$) because the requirement was
$a_{0}d_{0}-b_{0}c_{0}\neq 0$.\\
\newline \textbf{Case b.} If $a(x)$,
$b(x)$, $c(x)$ and $d(x)$ are arbitrary functions of $x$ then
\begin{eqnarray}
\tilde{y}^{\prime\prime}=a(x)y^{\prime\prime}+b(x)z^{\prime\prime}+a^{\prime\prime}(x)y+b^{\prime\prime}(x)z+2a^{\prime}(x)y^{\prime}+2b^{\prime}(x)z^{\prime},\nonumber\\
\tilde{z}^{\prime\prime}=c(x)y^{\prime\prime}+d(x)z^{\prime\prime}+c^{\prime\prime}(x)y+d^{\prime\prime}(x)z+2c^{\prime}(x)y^{\prime}+2d^{\prime}(x)z^{\prime}.\label{49}
\end{eqnarray}
Thus we obtain
\begin{eqnarray}
(ad-bc)y^{\prime\prime}=[(ad+bc)\tilde{d}_{11}+cd\tilde{d}_{12}-ab\tilde{d}_{21}-a^{\prime\prime}d+c^{\prime\prime}b]y+(2bd\tilde{d}_{11}+\nonumber\\
d^{2}\tilde{d}_{12}-b^{2}\tilde{d}_{21}-b^{\prime\prime}d+d^{\prime\prime}b)z-2d(a^{\prime}y^{\prime}+b^{\prime}z^{\prime})+2b(c^{\prime}y^{\prime}+d^{\prime}z^{\prime}),\label{50}\\
(ad-bc)z^{\prime\prime}=(2ac\tilde{d}_{11}+c^{2}\tilde{d}_{12}-a^{2}\tilde{d}_{21}-a^{\prime\prime}c+c^{\prime\prime}a)y+[(ad+bc)\tilde{d}_{11}+\nonumber\\
cd\tilde{d}_{12}-ab\tilde{d}_{21}-b^{\prime\prime}c+d^{\prime\prime}a]z-2c(a^{\prime}y^{\prime}+b^{\prime}z^{\prime})+2a(c^{\prime}y^{\prime}+d^{\prime}z^{\prime}).\label{51}
\end{eqnarray}
Comparing the coefficients as before and using the linear
independence of $\tilde{d}$'s we obtain
\begin{eqnarray}
a^{\prime}(x)=b^{\prime}(x)=c^{\prime}(x)=d^{\prime}(x)=0,\label{52}
\end{eqnarray}
which implies that it reduces to a system of the form ($\ref{47}$),
which leaves us again with the same result. Thus we have the
theorem.\\
\newline \textbf{Theorem 7.} \textit{The linear forms for
systems of two second-order ODEs obtainable by CSA are in general
inequivalent to those linear forms obtained by real symmetry
analysis.}

Before presenting some illustrative applications of the theory
developed we refine Theorem 6 by using Theorem 7 to make the
following remark.\\
\newline \textbf{Remark 2.} There are \emph{only}
$6$, $7$ or $15$-dimensional algebras for linearizable systems
obtainable by scalar complex linearizable ODEs, i.e. there are no
$5$ or $8$-dimensional Lie point symmetry algebras for such systems.

\section{Applications}
Consider a system of non-homogeneous geodesic-type differential
equations
\begin{eqnarray}
y^{\prime\prime}+y^{\prime 2}-z^{\prime 2}=\Omega_{1}(x,y,z,y^{\prime},z^{\prime}),\nonumber \\
z^{\prime\prime}+2y^{\prime}z^{\prime}=\Omega_{2}(x,y,z,y^{\prime},z^{\prime}).\label{53}
\end{eqnarray}
where $\Omega_{1}$ and $\Omega_{1}$ are linear functions of the
dependent variables and their derivatives. This system corresponds
to a complex scalar equation
\begin{eqnarray}
u^{\prime\prime}+u^{\prime 2}= \Omega(x,u,u^{\prime}),\label{54}
\end{eqnarray}
which is either transformable to the free particle equation or one
of the linear forms ($\ref{21}$)-($\ref{23}$), by means of the
complex transformations
\begin{eqnarray}
\chi=\chi(x), ~U(\chi)=e^{u}.\label{55}
\end{eqnarray}
Which are further transformable to the free particle equation by
utilizing another set of invertible complex point transformations.
Generally, the system ($\ref{53}$) is transformable to a system of
the free particle equations or a linear system of the form
\begin{eqnarray}
Y^{\prime\prime}=\widetilde{\Omega}_{1}(\chi, Y, Z, Y^{\prime}, Z^{\prime})-\widetilde{\Omega}_{2}(\chi, Y, Z, Y^{\prime}, Z^{\prime}),\nonumber \\
Z^{\prime\prime}=\widetilde{\Omega}_{2}(\chi, Y, Z, Y^{\prime},
Z^{\prime})+\widetilde{\Omega}_{1}(\chi, Y, Z, Y^{\prime},
Z^{\prime}).\label{56}
\end{eqnarray}
Here $\widetilde{\Omega}_{1}$ and $\widetilde{\Omega}_{2}$ are
linear functions of the dependent variables and their derivatives,
via an invertible change of variables obtainable from ($\ref{55}$).
The linear form ($\ref{56}$) can be mapped to a maximally symmetric
system if and only if there exist some invertible complex
transformations of the form ($\ref{55}$), otherwise these forms can
not be reduced further. This is the reason why we obtain three
equivalence classes namely with $6$, $7$ and $15$-dimensional
algebras for systems corresponding to linearizable complex equations
with only one equivalence class. We first consider an example of a
nonlinear system that admits a $15-$dimensional algebra which can be
mapped to the free particle system using ($\ref{55}$). Then we
consider four applications to nonlinear systems of quadratically
semi-linear ODEs transformable to ($\ref{56}$) via ($\ref{55}$) that
are not further
reducible to the free particle system.\\
\newline \textbf{1.} Consider ($\ref{53}$) with
\begin{eqnarray}
\Omega_{1}=-\frac{2}{x}y^{\prime}, \nonumber\\
\Omega_{2}=-\frac{2}{x}z^{\prime},\label{58}
\end{eqnarray}
it admits a $15$-dimensional algebra. The real linearizing
transformations
\begin{eqnarray}
\chi(x)=\frac{1}{x},~Y=e^{y}\cos(z),~Z=e^{y}\sin(z),\label{57}
\end{eqnarray}
obtainable from the complex transformations ($\ref{55}$) with
$U(\chi)=Y(\chi)+iZ(\chi)$, map the above nonlinear system to
$Y^{\prime\prime}=0$, $Z^{\prime\prime}=0$. Moreover, the solution
of ($\ref{58}$) corresponds to the solution of the corresponding
complex equation
\begin{eqnarray}
u^{\prime\prime}+u^{\prime 2}+\frac{2}{x}u^{\prime}=0.\label{59}
\end{eqnarray}
\newline
\textbf{2.} Now consider $\Omega_{1}$ and $\Omega_{2}$ to be linear
functions of the first derivatives $y^{\prime},~z^{\prime}$, i.e.,
system ($\ref{53}$) with
\begin{eqnarray}
\Omega_{1}=c_{1}y^{\prime}-c_{2}z^{\prime},\nonumber\\
\Omega_{2}=c_{2}y^{\prime}+c_{1}z^{\prime},\label{60}
\end{eqnarray}
which admits a $7$-dimensional algebra, provided both $c_{1}$ and
$c_{2}$, are not simultaneously zero. It is associated with the
complex equation
\begin{eqnarray}
u^{\prime\prime}+u^{\prime 2}-cu^{\prime}=0.
\end{eqnarray}
Using the transformations ($\ref{55}$) to generate the real
transformations
\begin{eqnarray}
\chi(x)= x,~Y=e^{y}\cos(z),~Z=e^{y}\sin(z),\label{61}
\end{eqnarray}
which map the nonlinear system to a linear system of the form
($\ref{24}$), i.e.,
\begin{eqnarray}
Y^{\prime\prime}=c_{1}Y^{\prime}-c_{2}Z^{\prime},\nonumber\\
Z^{\prime\prime}=c_{2}Y^{\prime}+c_{1}Z^{\prime},\label{62}
\end{eqnarray}
which also has a $7$-dimensional symmetry algebra and corresponds to
\begin{eqnarray}
U^{\prime\prime}-cU^{\prime}=0.\label{66}
\end{eqnarray}
All the linear second-order ODEs are transformable to the free
particle equation thus we can invertibly transform the above
equation to $\widetilde{U}^{\prime\prime}=0$, using
\begin{eqnarray}
(\chi(x), U)\rightarrow(\widetilde{\chi}=\alpha+\beta
e^{c\chi(x)},\widetilde{U}=U),
\end{eqnarray}
where $\alpha$, $\beta$ and $c$ are complex. But these complex
transformations can not generate real transformations to reduce the
corresponding system ($\ref{62}$) to a maximally symmetric system.\\
\newline \textbf{3.} A system with a $6-$dimensional Lie algebra is obtainable from ($\ref{53}$) by introducing a linear
function of $x$ in the above coefficients i.e.,
\begin{eqnarray}
\Omega_{1}=(1+x)(c_{1}y^{\prime}-c_{2}z^{\prime}),\nonumber\\
\Omega_{2}=(1+x)(c_{2}y^{\prime}+c_{1}z^{\prime}),\label{63}
\end{eqnarray}
in ($\ref{53}$), then the same transformations ($\ref{61}$) converts
the above system into a linear system
\begin{eqnarray}
Y^{\prime\prime}=(1+\chi) \left (c_{1}Y^{\prime}-c_{2}Z^{\prime}
\right ),\nonumber\\
 Z^{\prime\prime}=(1+\chi) \left
(c_{2}Y^{\prime}+c_{1}Z^{\prime} \right ),\label{64}
\end{eqnarray}
where both systems ($\ref{63}$) and ($\ref{64}$) are in agreement on
the dimensions (i.e. six) of their symmetry algebras. Again, the
above system is a special case of the linear system
($\ref{24}$).\\
\newline \textbf{4.} If we choose $\Omega_{1}=c_{1},~
\Omega_{2}=c_{2}$, where $c_{i}$ $(i=1,2)$ are non-zero constants,
then under the same real transformations ($\ref{61}$), the nonlinear
system ($\ref{53}$) takes the form
\begin{eqnarray}
Y^{\prime\prime}=c_{1}Y-c_{2}Z,\nonumber\\
Z^{\prime\prime}=c_{2}Y+c_{1}Z.\label{65}
\end{eqnarray}

%%%%%%%%%%%%%%%%%%%%%%%%%%%%%%%%%%%%%%%%%%%%%%%%%%%%%%%%%%%%%%%%%%%%%%%%%%%%%%%%%%%%%%%%%%%%%%%%%%%%%%%%%%%%%%%%%%%%%%%
%%%%%%%%%%%%%%%%%%%%%%%%%%%%%%%%%%%%%%%%%%%%%%%%%%%%%%%%%%%%%%%%%%%%%%%%%%%%%%%%%%%%%%%%%%%%%%%%%%%%%%%%%%%%%%%%%%%%%%%
\section{Conclusion}
The classification of linearizable systems of two second-order ODEs
was obtained by using the equivalence properties of systems of two
linear second-order ODEs \cite{sb}. The ``optimal canonical form" of
the corresponding linear systems of two second-order ODEs, to which
a linearizable system could be mapped, is crucial. This canonical
form used invertible transformations, the invertibility of these
mappings insuring that the symmetry structure is preserved. That
optimal canonical form of the linear systems of two second-order
ODEs led to five linearizable classes with respect to Lie point
symmetry algebras with dimensions $5$, $6$, $7$, $8$ and $15$.

Systems of two second-order ODEs appearing in CSA correspond to some
scalar complex second-order ODE. We proved the existence of a
reduced optimal canonical form for such linear systems of two ODEs.
This reduced canonical form provided three equivalence classes,
namely with $6$, $7$ or $15$-dimensional point symmetry algebras.
Two cases are eliminated in the theory of complex symmetries: those
of $5$ and $8$-dimensional algebras. The systems corresponding to a
complex linearized scalar ODE involve one parameter which can only
cover {\it three} possibilities; (a) it is zero; (b) it is a
non-zero constant; and (c) it is a non-constant function. The non
existence of $5$ and $8$ dimensional algebras for the linear forms
appearing due to CSA has been proved by showing that these forms are
not equivalent to those provided by the real symmetry approach for
systems \cite{sb} with $5$ and $8$ generators.

Work is in progress \cite{saf3} to find complex methods of solving a
class of 2-dimensional \emph{nonlinearizable} systems of
second-order ODEs. It is also obtainable from the linearizable
scalar complex second-order ODEs, which are transformable to the
free particle equation via an invertible change of the dependent and
independent variables of the form
\begin{eqnarray}
\chi=\chi(x,u), ~U(\chi)=U(x,u).
\end{eqnarray}
Notice that these transformations are different from ($\ref{55}$).
The real transformations corresponding to the complex
transformations above cannot be used to linearize the real system.
But the linearizability of the complex scalar equations can be used
to provide solutions for the corresponding systems.\\
%%%%%%%%%%%%%%%%%%%%%%%%%%%%%%%%%%%%%%%%%%%%%%%%%%%%%%%%%%%%%%%%%%%%%%%%%%%%%%%%%%%%%%%%%%%%%%%%%%%%%%%%%%%%%%%%%%%%%%%
%%%%%%%%%%%%%%%%%%%%%%%%%%%%%%%%%%%%%%%%%%%%%%%%%%%%%%%%%%%%%%%%%%%%%%%%%%%%%%%%%%%%%%%%%%%%%%%%%%%%%%%%%%%%%%%%%%%%%%%
\newline\textbf{Acknowledgements}\newline The authors are grateful to
Fazal Mahomed for useful comments and discussion on this work. MS is
most grateful to NUST for providing financial support.

%%%%%%%%%%%%%%%%%%%%%%%%%%%%%%%%%%%%%%%%%%%%%%%%%%%%%%%%%%%%%%%%%%%%%%%%%%%%%%%%%%%%%%%%%%%%%%%%%%%%%%%%%%%%%%%%%%%%%%%%%
%SECTION
\bc{\bf REFERENCES}\ec \vspace{-1.8cm}

\renewcommand{\refname}{}


\begin{thebibliography}{99}
\bibitem{saj} S. Ali, F.M. Mahomed and A. Qadir, Criteria for systems of
two second-order differential equations by complex methods, {\it
Nonlinear Dynamics} (to appear).

\bibitem{saj1} S. Ali, F.M. Mahomed and A. Qadir, Complex
Lie symmetries for Scalar Second-order Ordinary Differential
Equations, {\it Nonlinear Analysis: Real World Applications}
\textbf{10} (2009) 3335-3344.

\bibitem{saf} S. Ali, F.M. Mahomed and M. Safdar, Complete
classification for systems of two second-order ordinary differential
equations by using complex symmetry analysis, work in progress.

\bibitem{saf3} S. Ali, A. Qadir and M. Safdar, Symmetry solutions of two dimensional systems not solvable by symmetry methods, work in
progress.


\bibitem{aa} A.V. Aminova, N.A.M. Aminov, Projective geometry of
systems of differential equations: general conceptions, {\it Tensor
N S} {\bf 62} (2000), 65-86.

\bibitem{gor} V.M. Gorringe and P. G.L. Leach, Lie point
symmetries for systems of $2$nd order linear ordinary differential
equations, {\it Quaestiones Mathematicae} \textbf{11}(1) (1988)
95-117.

\bibitem{ibr} N.H. Ibragimov, {\it Elementary Lie Group Analysis and
Ordinary Differential Equations}, Wiley Chichester (1999).

\bibitem{lie1} S. Lie, {\it Differential Equations}, Chelsea, New
York (1967).

\bibitem{lie} S. Lie, Klassifikation und Integration von
gew\"onlichen Differentialgleichungenzwischen $x$, $y$, die eine
Gruppe von Transformationen Gestaten, {\it  Arch. Math.}  {\bf VIII,
IX} (1883), 187.

\bibitem{lie2} S. Lie, {\it Lectures on Differential Equations with
Known Infinitesimal Transformations}, Teubner: Leipzig (1891) (in
German, Lie's lectures by G. Sheffers).

\bibitem{lie3} S. Lie, Theorie der Transformationsgruppen,
{\it Math. Ann.}, {\bf 16} (1880), 441.

\bibitem{sur} F.M. Mahomed, Symmetry group classification of
ordinary differential equations: Survey of some results, {\it Math.
Meth. Appl. Sci.} \textbf{30} (2007) 1995-2012.

\bibitem{lea} F.M. Mahomed and P.G.L. Leach, Symmetry Lie
algebra of $n$th order ordinary differential equations, {\it J.
Math. Anal. Appl.} \textbf{151} (1990) 80-107.

\bibitem{mq1} F.M. Mahomed and A. Qadir, Linearization criteria for a
system of second-order quadratically semi-linear ordinary
differential equations, {\it Nonlinear Dynamics} {\bf 48} (2007)
417-422.

\bibitem{mq2} F.M. Mahomed and A. Qadir, Invariant
linearization criteria for systems of cubically semi-linear
second-order ordinary differential equations, {\it J. Nonlinear
Math. Phys.} \textbf{16} (2009) 1-16.

\bibitem{jm} J. Merker, Characterization of the Newtonian free
particle system in $m\geq2$ unknown dependent variables, {\it Acta
Applicandae Mathematicae} \textbf{92} (2006) 125-207.

\bibitem{saf2} M. Safdar, A. Qadir and S. Ali, Inequivalence of
classes of linearizable systems of cubically semi-linear ordinary
differential equations obtained by real and complex symmetry
analysis, {\it Math. Comp. Appl.} (to appear).

\bibitem{hs} H. Stephani, {\it Differential Equations: Their
Solutions Using Symmetries}, Cambridge (1996).

\bibitem{sb} C. Wafo Soh and F.M. Mahomed, Symmetry breaking
for a system of two linear second-order ordinary differential
equations, {\it Nonlinear Dynamics} \textbf{22} (2000) 121-133.

\bibitem{wf} C. Wafo Soh and F.M. Mahomed, Linearization
criteria for a system of second-order ordinary differential
equations, {\it Int. J. of NonLinear Mech.} {\bf36} (2001) 671-677.


\end{thebibliography}
\end{document}